\input amstex
\input table
\documentstyle{amsppt}


\ifx\undefined\headlineheight\else\pageheight{47pc}\vcorrection{-1pc}\fi
 
\magnification=\magstep1
\pagewidth{6.5truein}
\pageheight{8.9truein}
\NoBlackBoxes

\def\pf{\qed}

\topmatter
\title
	Some Results on Combinators in the System TRC
\endtitle
\author
Thomas Jech
\endauthor
\affil 
The Pennsylvania State University
\endaffil
\thanks
Supported in part by the National Science Foundation grant
DMS--9401275.  A computer equipment provided by the NSF SCREMS Grant
DMT-9628630.  Some of these results were presented at the Special
Session on Computer Proofs in Logic at the AMS meeting in Atlanta in
October 1997.
\endthanks
\address
Department of Mathematics, The Pennsylvania State University,
218 McAllister Bldg., University Park, PA 16802, U.S.A. 
\endaddress
\email jech\@math.psu.edu \endemail
\abstract
We investigate the system TRC of untyped illative combinatory
logic that is equiconsistent with New Foundations.  We prove
that various unstratified combinators do not exist in TRC.
\endabstract
\endtopmatter

\document
\baselineskip=18truept
\subhead
{Introduction}
\endsubhead
We prove some results in the axiomatic system TRC introduced in
\cite3.  The system TRC (for `type-respecting combinators') is an
untyped system of combinatory logic, in the sense of \cite1, \cite2.  TRC is a
first order theory of functions (combinators) with equality and is
illative, i.e. capable of expressing notions of propositional logic.
Moreover, it is combinatorially complete for stratified combinators.
The main interest of TRC is that it is equiconsistent with the theory
NF \cite6, Quine's `New Foundations'.  As the consistency of NF
remains an open problem, so does the consistency of TRC.

The objects of study of a combinatory logic are {\it combinators}.  We
denote $xy$ the application of the combinator $x$ to the combinator
$y$, and adopt the convention that $xyz = (xy)z$.

The language of TRC has (in addition to equality and the binary
function $xy$) constants $Abst$, $Eq$, $p_1$ and $p_2$, and functions
$k(x)$ and $\langle x,y\rangle$.  The axioms of TRC are the following:
\roster 
\item"{I. }" $k(x) y = x$.
\item"{II. }" $p_i \langle x_1,x_2\rangle = x_i$ for $i = 1,2$.
\item"{III. }" $\langle p_1 x, p_2 x\rangle = x$.
\item"{IV. }" $\langle x,y\rangle z = \langle xy,xz\rangle$
\item"{V. }" $Abst\,x\,y\,z = x\, k(y)(yz)$.
\item"{VI. }"$ Eq \langle x,y\rangle = p_1$ if $x = y$; $Eq 
	\langle x,y\rangle = p_2$ if $x \neq y$.
\item"{VII. }" If for all $z$, $xz = yz$, then $x = y$.
\item"{VIII. }" $p_1 \neq p_2$.
\endroster
 Axiom I postulates the existence of constant functions.  Axioms
 II--IV describe the pairing function $\langle x,y\rangle$ and the
 projections $p_1$ and $p_2$.  $Abst$ is the abstraction combinator and
 $Eq$ is the characteristic function of equality.  Axiom VII is the
 axiom of extensionality.

Let $I = \langle p_1,p_2\rangle$; from III and IV it follows 
that $I$ is the identity function $Ix = x$.

Classical combinatory logic \cite{2,3} employs combinators
$I$, $K$ and $S$, where
$$
	Ix = x\quad\quad Kxy = x\quad\quad Sxyz = xz(yz).
$$
It has a powerful abstraction property: for every term $t$ and a
variable $x$, there is a term $\lambda xt$ in which $x$ does not
occur, such that for every term $s$,
$$
	(\lambda xt)s = t[s/x].
$$
This guarantees, among others, the existence of a fixed point
for every combinator, and implies that simple notions of
propositional logic cannot be represented by combinators.  
Suppose that $Neg$ is the negation combinator, and consider
$u = \lambda x(Neg(xx))$.  Then $uu = Neg(uu)$.

The theory TRC is an illative theory, in the sense that it can encode
notions of propositional logic.  It also has an abstraction property
(Theorem 1 of \cite3).  The term $\lambda xt$ can be constructed for
every $t$ in which $x$ occurs with no type other than $0$.  (For
details about typing see \cite3.)  It follows that TRC proves the
existence of all {\it stratified} combinators.  Examples of stratified
combinators are $x(yx)$, $xy(yz)$, $y(xyz)$: in $y(xyz)$, $z$ has type
$0$, $y$ has type $1$ and $x$ has type $2$.  (In fact, $Abst\,Abst\,Ixy =
x(yx)$, $Abst(AbstAbst)xyz = xy(yz)$, and $AbstAbst\, xyz = y(xyz)$).

We will show in Section 3 that (with the exception of $I$) the
standard combinators used in classical combinatory logic do not exist
in TRC.  We shall give many examples of unstratified combinators whose
existence contradicts the axioms of TRC.

In searching for proofs of the various results in TRC, we used a
computer extensively and used the automated theorem prover OTTER
\cite{5}.

\subhead
2. Some Basic Facts on TRC
\endsubhead

In this sectioin we derive some simple equalities from the axioms
of TRC, and use a self-reference argument to obtain some simple
negative results.  First we state some properties of the abstraction
combinator  (see also \cite{4}):

\proclaim
{Theorem 2.1}
\roster
\runinitem"{(a) }" $Abst(Abst(Abst\;x)) = Abst\, x$.
\item"{(b) }"  $Abst(Abst \,k(x)) = k(x)$.
\item"{(c) }"  $Abst \;k(k(x)) = k(k(x))$.
\item"{(d) }"  $Abst\; k(x) yz = x(yz)$.
\item"{(e) }"  $Abst\; k(x)k(y) = k(xy)$.
\endroster
\endproclaim

\demo
{Proof}  The equalities are obtained by an application of the axioms
defining $Abst$ and $k(x)$ and the axiom of extensionality; e.g. to prove
(a), we evaluate the term $Abst(Abst(Abst\;x))yz$ and compare it
with $Abst\;x\,y\,z$. \newline
\line{\hfil \pf}
\enddemo

The next theorem gives some properties of the pairing function
and the projections:

\proclaim
{Theorem 2.2}
\roster
\runinitem"{(a) }" $\langle k(x),y\rangle z = \langle x,yz\rangle,
\langle x,k(y)\rangle z = \langle  xz,y\rangle$.
\item"{(b) }"  $k(\langle x,y\rangle) = \langle k(x),k(y)\rangle,\;
	k(p_i x) = p_i k(x)$ for $i = 1,2$.
\item"{(c) }"  $p_i(xy) = p_i xy$ for $i = 1,2$.
\endroster
\endproclaim

\demo
{Proof}  \roster
\runinitem"{(a) }" From Axiom IV.
\item"{(b) }"  Calculate $\langle k(x),k(y)\rangle z$ and $p_i k(x)z$
	and use extensionality.
\item"{(c) }"  Let $x = \langle u,v\rangle$ and use Axiom IV.  \pf
\endroster
\enddemo

Next we state some more properties of the combinator $Abst$:

\proclaim
{Theorem 2.3} \roster
\runinitem"{(a) }" $Abst \langle x,y\rangle = \langle Abst 
	\;x,Abst\;y\rangle$.
\item"{(b) }"  $Abst\; p_i = k(p_i)$ and $Abst\;k(p_i) = p_i$, for 
	$i = 1,2$.
\item"{(c) }"  $Abst\; I = k(I)$ and $Abst\;k(I) = I$.
\endroster
\endproclaim

\demo
{Proof}  \roster
\runinitem"{(a) }" Using Axiom IV, show that $Abst \langle x,y\rangle 
	uv = \langle Abst\;x,Abst\;y\rangle uv$.
\item"{(b) }"  $Abst \;p_i xy = p_i k(y)(xy) = k(p_i y)(xy) = p_i y$ 
	by Theorem 2.2b, and $Abst \;k(p_i)xy = p_i(xy) = p_i xy$ by 
	Theorems 2.1d and 2.2c.
\item"{(c) }"  $Abst\,Ixy = y = k(I) xy$, by Axioms V and I,
	and $Abst\; k(I) xy = I(xy) = Ixy$ by Theorem 2.1d.  \pf
\endroster
\enddemo

We shall now turn to negative results.  In Section 3 we 
shall present a number of combinators that do not exist in
TRC.  Each proof will use one of the following basic negative
results that use self-reference:

\proclaim
{Theorem 2.4}  For every $x$,
\roster
\item"{(a) }"  $Eq \langle x,p_2\rangle \neq x$.
\item"{(b) }"  $Eq \langle k(x),k(p_2)\rangle \neq x$.
\item"{(c) }"  $\langle Eq\,x,p_2\rangle \neq x$.
\endroster
\endproclaim

\demo
{Proof} (a) $Eq\langle x,y\rangle$ is either $p_1$ or $p_2$, and
$Eq\langle p_1,p_2\rangle = p_2$ while $Eq\langle p_2,p_2\rangle =
p_1$.

The proof of (b) and (c) is similar.  \pf
\enddemo

It follows from the discussion on classical combinatory logic in 
Section 1 that not every combinator in TRC has a fixed point.
Theorem 2.4 gives an explicit example, $\langle Eq, k(p_2)\rangle
x \neq x$:

\proclaim
{Corollary 2.5}  The combinator $\langle Eq, k(p_2)\rangle$
does not have a fixed point.
\endproclaim

A standard fact of combinatory logic (cf. \cite1, \cite2) states that
if $M$ is the combinator $Mx = xx$ then for every $u$, the
composition of $u$ and $M$ is a fixed point of $u$.  As the 
abstraction theorem for TRC in \cite3 provides for composition
of combinators, it follows that $M$ does not exist in TRC.
Here we give a direct proof:

\proclaim
{Theorem 2.6}
There is no $M$ such that $Mx = xx$.
\endproclaim

\demo
{Proof}  Let $t = Abst \,k(Eq) \langle M, k(p_2)\rangle$ and  let
$s = tt$.  Then (using Theorems 2.1.d and 2.2a)
$$
\aligned
	s = tt & = Abst\,k (Eq) \langle M,k(p_2)\rangle t\\
	& = Eq(\langle M,k(p_2)\rangle t)   \\
	& = Eq \langle Mt,p_2\rangle   \\
	& = Eq \langle tt, p_2\rangle \\
	& = Eq \langle s,p_2\rangle,
\endaligned
$$
contradicting Theorem 2.4a.  \pf
\enddemo

A similar argument, using Theorem 2.4b, yields the following:

\proclaim
{Theorem 2.7}  There is no $K_1$ such that $K_1 x = k(xx)$.
\endproclaim

\demo
{Proof}  Let $t = Abst\,k(Eq) \langle K_1,k(k(p_2))\rangle$, and
$s = tt$.  Then (by Theorems 2.1d and 2.2a)
$$
\aligned
	s = tt & = Abst\, k(Eq) \langle K_1,k(k(p_2))\rangle t \\
	& = Eq(\langle K_1,k(k(p_2))\rangle t)   \\
	& = Eq \langle K_1 t,k(p_2)\rangle   \\
	& = Eq \langle k(s),k(p_2)\rangle,
\endaligned
$$
contradicting Theorem 2.4b.  \pf
\enddemo

An immediate consequence of Theorem 2.7 is that neither $k(x)$ nor
$\langle x,y\rangle$ can be replaced in TRC by a combinator (see
also \cite4).

\proclaim
{Theorem 2.8} \roster
\runinitem"{(a) }" There is no $K$ such that $Kx = k(x)$.
\item"{(b) }"  There is no $p$ such that $pxy = \langle x,y\rangle$.
\endroster
\endproclaim

\demo
{Proof}\roster
\runinitem"{(a) }"  Given such $K$, let $K_1 = AbstAbst \;K$.
Then (by Theorem 2.1d)
$$
\aligned
	K_1 xy & = Abst Abst\;Kxy  \\
	& = Abst \;K(x)(Kx)y   \\
	& = x(Kxy)   \\
	& = xx
\endaligned
$$
and so $K_1 x = k(xx)$, contradicting Theorem 2.7.
\item"{(b) }"  Given $p$, let $K = p_1 p$, and then (by Theorem 2.2c)
$$
\aligned
	Kxy & = p_1 pxy  \\
	& = p_1(px)y   \\
	& = p_1(pxy)   \\\
	& = p_1\langle x,y\rangle \\
	& = x,
\endaligned
$$
contradicting (a).   \pf
\endroster
\enddemo

We conclude this section with the following result that we use
in Section 3.

\proclaim
{Theorem 2.9}  \roster 
\runinitem"{(a) }" There is no $u$ such that $ux = xk(x)$.
\item"{(b) }"  There is no $u$ such that $uk(x) = xk(x)$.
\endroster
\endproclaim

\demo
{Proof}  \roster
\runinitem"{(a) }" Given $u$, let $M = Abst(Abst\;u)I$, and
then
$$
\aligned
	Mx & = Abst(Abst \;u)I\,x   \\
	& = Abst\,u\,k(x)\,x   \\
	& = u\,k(x)\,x   \\
	& = k(x) k(k(x))x   \\
	& = x\,x,
\endaligned
$$
contradicting Theorem 2.6.
\item"{(b) }"  Given $u$, let $t = Abst\;k(Eq)\langle u,k(p_2)\rangle$
and $s = uk(t)$.  Then we have (by Theorems 2.1.d and 2.2a)
$$
\aligned 
	s = u\,k(t) & = Abst\,k(Eq) \langle u,k(p_2)\rangle k(t) \\
	& = Eq(\langle u,k(p_2)\rangle k(t))   \\
	& = Eq \langle u\,k(t),p_2\rangle   \\
	& = Eq \langle s,p_2 \rangle,
\endaligned
$$
contradicting Theorem 2.4a.  \pf
\endroster
\enddemo

\subhead
3. Nonexistence of Various Combinators
\endsubhead

We will show that many standard classical combinators do not
exist in TCR.  Let us consider the following combinators; 
none of them is stratified.  We use the list presented in \cite7,
with several additions.

$$
\begintable
\begintableformat
\center "  {\qquad\quad\center}  "  \center " {\qquad\quad\center}  " \center 
\endtableformat
\br{\:} $B xyz = x(yz)$ " {} " $Lxy = x(yy)$ " {} " $Q_3xyz = z(xy)$ \er{}
\br{\:} $C xyz = xzy$ " {} " $L_1xy = y(xx)$ " {} " $Rxyz = yzx$ \er{}
\br{\:} $D xyzw = xy(zw)$ " {} " $Mx = xx$ " {} " $Sxyz = xz(yz)$ \er{}
\br{\:} $F xyz = zyx$ " {} " $M_1 x = xxx$ " {} " $Txy = yx$ \er{}
\br{\:} $G xyzw = xw(yz)$ " {} " $M_2 x = x(xx)$ " {} " $Uxy = y(xxy)$ \er{}
\br{\:} $H xyz = xyzy$ " {} " $Oxy = y(xy)$ " {} " $Vxyz = zxy$ \er{}
\br{\:} $H_1 xy = xyx$ " {} " $O_1 xy = x(yx)$ " {} " $Wxy=xyy$ \er{}
\br{\:} $J xyzw = xy(xwz)$ " {} " $O_2 xy = y(yx)$ " {} " $W_1xy = yxx$ \er{}
\br{\:} $K xy = x$ " {} " $Q xyz = y(xz)$ " {} " $W_2xy = yxy$ \er{}
\br{\:} $K_1 xy = xx$ " {} " $Q_1 xyz = x(zy)$ " {} " $W_3xy = yyx$ \er{}
\endtable
$$

Below we prove that none of these combinators exist in TRC.

\subhead
(3.1)
\endsubhead
$K_1$, $K$, $M$ and $J$:

Theorems 2.6, 2.7 and 2.8 show that $K_1$, $K$ and $M$ do not exist.
As for $J$, it is well known in combinatory logic (cf. \cite1) that
$\{I,J\}$ is combinatorially complete, and so $J$ cannot exist 
in TRC.

\subhead
(3.2)
\endsubhead
$L$, $O$, $U$ and $W$:
$$
	M = LI = OI = UI = WI
$$

\subhead
(3.3)
\endsubhead
$O_2$ and $M_2$:
$$
\aligned
	M = Abst(O_2 I) I & = Abst\; M_2 I:  \\
	Abst (O_2 I)Ix & = O_2 I\,k(x)x = k(x) (k(x)I) x = xx  \\
	Abst\;M_2 Ix & = M_2 k(x) x = k(x)(k(x)k(x))x = xx
\endaligned
$$

\subhead
(3.4)
\endsubhead
$S$ and $O_1$:

$O = SI$ and $S = Abst \circ O_1$ (where $a \circ b$ is the
composition, defined in TRC by $a \circ b = Abst\;k(a)(Abst\;k(b)I)$):
$$
      Sxyz = Abst(O_1x)yz = O_1 x\,k(z)(yz) = x(k(z)x)(yz) = xz(yz)
$$

\subhead
(3.5)
\endsubhead
$T$, $C$, $G$, $Q_1$ and $Q_3$:
$$
\gathered
	K = Abst\;T\,k(I) \text{ and } T = CI = GII = Q_1 I = Q_3 I:  \\
	Kx = Abst\;T\, k(I)x = T\;k(x)(k(I)x) = T\;k(x) I = Ik(x)
	= k(x)
\endgathered
$$

\subhead
(3.6)
\endsubhead
$B$ and $D$:
$$
\aligned
	K & = Abst BI, \; B = DI: \\
	K xy & = Abst\;B\,I\,xy = B\,k(x)(Ix)y = B\,k(x) xy =  \\
	& = k(x) (xy) = x
\endaligned
$$

\subhead
(3.7)
\endsubhead
$R$:

$K = R\,k(I) p_1 \langle R,u\rangle R$, where $u$ is arbitrary:
$$
\aligned
	K xy & = R\,k(I) p_1 \langle R,u\rangle R xy  \\
	& = p_1 \langle R,u\rangle k(I) R\,xy  \\
	& = R\,k(I)\, R\,xy  \\
	& = R\,x\,k(I)y    \\
	& = k(I) yx  \\
	& = Ix = x
\endaligned
$$

\subhead
(3.8)
\endsubhead
$V$:
$$
	K = Abst (Abst\;V\;Abst)k(k(Abst)):
$$
using Theorem 2.1.b, we have
$$
\aligned
	Kx & = Abst (Abst\;V\;Abst) k(k(Abst))x  \\
	& = Abst\;V\; Abst\; k(x) k(Abst)   \\
	& = V\;k(k(x)) (Abst \;k(x)) k(Abst)   \\
	& = k(Abst) k(k(x)) (Abst\; k(x))   \\
	& = Abst (Abst\;k(x))  \\
	& = k(x)
\endaligned
$$

\subhead
(3.9)
\endsubhead
$Q$:
$$
\aligned
	K_1 & = Abst\;Q\,I:  \\
	K_1 xy & = Abst\; Q\,I\,xy = Q\;k(x)xy = x(k(x)y) = xx
\endaligned
$$

\subhead
(3.10)
\endsubhead
$H_1$, $H$, $M_1$ and $W_2$:
$$
	M_1 k(x) = H_1 H_1 k(x) = W_2 W_2 k(x) = xk(x),
$$
contradicting Theorem 2.9, and $H_1 = HI$.

\subhead
(3.11)
\endsubhead
$F$ and $W_1$:

Let $u = Abst(Fz) Abst \;k(x)$ (where $z$ is arbitrary) and
$v = Abst\; W_1\;Abst$.  Then $u\,k(x) = v\,k(x) = x\,k(x)$, 
contradicting Theorem 2.9:  using Theorem 2.1d, we have
$$
\aligned
	uk(x) & = Abst(Fz)Abst \;k(x)  \\
	& = Fz\; k(k(x)) (Abst\;k(x))  \\
	& = Abst\;k(x)\;k(k(x))\,z  \\
	& = x(k(k(x)) z)  \\
	& = x\,k(x)
\endaligned
$$
and
$$
\aligned
	vk(x) & = Abst \;W_1\;Abst\;k(x)  \\
	& = W_1\;k(k(x))\,(Abst\;k(x))  \\	
	& = Abst\;k(x)\;k(k(x))\;k(k(x))  \\
	& = x\,(k(k(x))\;k(k(x)))   \\
	& = x\,k(x).
\endaligned
$$

\subhead
(3.12)
\endsubhead
$L_1$:

Let $a = k(\langle Eq,k(p_2)\rangle)$.  Then for all $x$,
$$
\aligned
	L_1 a(L_1 x) & = L_1 x(aa)  \\
	& = aa(xx)   \\
	& = k(\langle Eq,k(p_2)\rangle) a(xx)   \\
	& = \langle Eq,k(p_2)\rangle(xx)  \\
	& = \langle Eq (xx), p_2\rangle,
\endaligned
$$
which, by Theorem 2.4c, is not equal to $xx$.  

Now let $b = Abst\,k(L_1 a) L_1$.  By Theorem 2.1d we have
$$
	bb = Abst\;k(L_1 a) L_1 b = L_1 a(L_1 b),
$$
a contradiction.

\subhead
(3.13)
\endsubhead
$W_3$:

Let $a = k(\langle Eq, k(p_2)\rangle)$.  Then for all $x$,
$$
\aligned
	W_3 x\, a & = a\,a\,x   \\
	& = \langle Eq, k(p_2)\rangle x  \\
	& = \langle Eq\; x, p_2\rangle, \\
\endaligned
$$
which, by Theorem 2.4c, is not equal to $x$.

Now let $b = Abst\;k(W_3)(W_3 a)$.  By Theorem 2.1d we have
$$
	W_3 ab = bba = Abst\;k(W_3)(W_3 a)ba = W_3 (W_3 ab)a.
$$
Thus if above we let $x = W_3 ab$, we get
$$
	W_3 (W_3 ab) a \neq W_3 ab,
$$
a contradiction.

\enddocument \end